\documentclass{amsart}
\theoremstyle{plain}
  \newtheorem{thm}{Theorem}[section]
  \newtheorem{lem}[thm]{Lemma}
  
  \newtheorem{prop}[thm]{Proposition}
  \newtheorem{conj}[thm]{Conjecture}

\theoremstyle{definition}

\theoremstyle{remark}
  \newtheorem{rem}[thm]{Remark}
  \newtheorem*{ack}{Acknowledgments}
\newcommand{\cR}{\check{R}}

\newcommand{\tR}{\tilde{R}}
\newcommand{\Z}{\mathbb{Z}}
\newcommand{\C}{\mathbb{C}}
\newcommand{\bH}{\mathbb{H}}
\newcommand{\qbinom}[2]{\genfrac{[}{]}{0pt}{}{#1}{#2}}
\newcommand{\res}{\operatorname{res}}
\newcommand{\Vol}{\operatorname{Vol}}

\newcommand{\End}{\operatorname{End}}
\renewcommand{\Sp}{\mathrm{Sp}}
\numberwithin{equation}{section}
\date{\today}
\allowdisplaybreaks
\begin{document}
\title[The colored Jones polynomials and the simplicial volume of a knot]
{The colored Jones polynomials \\ and \\ the simplicial volume of a knot}
\author{Hitoshi Murakami}
\address{
  Department of Mathematics,
  School of Science and Engineering,
  Waseda University,
  Ohkubo, Shinjuku, Tokyo 169-8555, Japan
  \\
  and
  Mittag-Leffler Institute,
  Aurav{\"a}gen 17,
  S-182 62, Djursholm,
  Sweden
}
\email{hitoshi@uguisu.co.jp}
\author{Jun Murakami}
\address{
  Department of Mathematics,
  Graduate School of Science,
  Osaka University,
  Machi\-kane\-yama-cho 1-1,
  Toyonaka, Osaka 560-0043, Japan
}
\email{jun@math.sci.osaka-u.ac.jp}
\thanks{
This research is supported in part by Sumitomo Foundation and Grand-in-Aid for
Scientific Research, The Ministry of Education, Science, Sports and Culture.
}
\begin{abstract}
We show that the set of colored Jones polynomials and the set of
generalized Alexander polynomials defined by Akutsu, Deguchi and Ohtsuki
intersect non-trivially.
Moreover it is shown that the intersection is (at least includes) the set of
Kashaev's quantum dilogarithm invariants for links.
Therefore Kashaev's conjecture can be restated as follows:
The colored Jones polynomials determine the hyperbolic volume for a hyperbolic
knot.
Modifying this, we propose a stronger conjecture:
The colored Jones polynomials determine the simplicial volume for any knot.
If our conjecture is true, then we can prove that a knot is trivial if and only
if all of its Vassiliev invariants are trivial.
\end{abstract}
\maketitle
In \cite{Kashaev:MODPLA95}, R.M.~Kashaev defined a family of complex valued link
invariants indexed by integers $N\ge2$ using the quantum dilogarithm.
Later he calculated the asymptotic behavior of his invariant and observed
that for the three simplest hyperbolic knots it grows as
$\exp(\Vol(K)N/2\pi)$ when $N$ goes to the infinity, where $\Vol(K)$ is the
hyperbolic volume of the complement of a knot $K$ \cite{Kashaev:LETMP97}.
This amazing result and his conjecture that the same also holds for any
hyperbolic knot have been almost ignored by mathematicians since his definition
of the invariant is too complicated (though it uses only elementary tools).
\par
The aim of this paper is to reveal his mysterious definition and to show that
his invariant is nothing but a specialization of the colored Jones polynomial.
The colored Jones polynomial is defined for colored links
(each component is decorated with an irreducible representation of the Lie
algebra $sl(2,\C)$).
The original Jones polynomial corresponds to the case that all the colors are
identical to the $2$-dimensional fundamental representation.
We show that Kashaev's invariant with parameter $N$ coincides with the colored
Jones polynomial in a certain normalization with every color the $N$-dimensional
representation, evaluated at the primitive $N$-th root of unity.
(We have to normalize the colored Jones polynomial so that the value for
the trivial knot is one, for otherwise it always vanishes).
\par
On the other hand there are other colored polynomial invariants, the generalized
multivariable Alexander polynomial defined by Y.~Akutsu, T.~Deguchi and
T.~Ohtsuki \cite{Akutsu/Deguchi/Ohtsuki:JKNOT92}.
They used the same Lie algebra $sl(2,\C)$ but a different hierarchy of
representations.
Their invariants are parameterized by $c+1$ parameters; an integer $N$ and
complex numbers $p_i$ ($i=1,2,\dots,c$) decorating the components,
where $c$ is the number of components of the link.
In the case where $N=2$, their invariant coincides with the multivariable
Alexander polynomial and their definition is the same as the second
authors' \cite{Jun:PACJM93}.
Using the Akutsu--Deguchi--Ohtsuki invariants we have another coincidence.
We will show that if all the colors are $(N-1)/2$ then the generalized Alexander
polynomial is the same as Kashaev's invariant since it coincides with the
specialization of the colored Jones polynomial as stated above.
Therefore the set of colored Jones polynomials and the set of generalized
Alexander polynomials of Akutsu--Deguchi--Ohtsuki intersect at
Kashaev's invariants.
\par
The paper is organized as follows.
In the first section we recall the definition of the link invariant defined
by Yang--Baxter operators.
In \S\ref{sec:Alexander} we show that the Akutsu--Deguchi--Ohtsuki invariant
coincides with the colored Jones polynomial when the colors are all
$(N-1)/2$ by showing that their representation becomes the usual representation
corresponding to the irreducible $N$-dimensional representation of $sl(2,\C)$.
In \S\ref{sec:Jones} we show that if we transform the $\cR$-matrix
used in the colored Jones polynomial by a Vandermonde matrix then it has
a form very similar to Kashaev's $\cR$-matrix.
In fact it is proved in \S\ref{sec:Kashaev} that these two
$\cR$-matrix differ only by a constant.
We also confirm the well-definedness of Kashaev's invariant by using this fact.
\par
In the final section we propose our "dr{\"o}m i Djursholm".
We use M.~Gromov's simplicial volume for a knot to generalize Kashaev's
conjecture.
Observing that the simplicial volume is additive and unchanged by mutation,
we conjecture that Kashaev's invariants
($=$ specializations of the colored Jones polynomial $=$ specializations of the Akutsu--Deguchi--Ohtsuki invariant) determine the simplicial volume.
If our dream comes true, then we can show that a knot is trivial if and only
if all of its Vassiliev invariants are trivial.
\par
\begin{ack}
Most of this work was done during the authors were visiting the Mittag-Leffler
Institute, Djursholm.
We thank the Institute and the staff for their hospitality.
We are also grateful to D.~Bar-Natan, E.~Date, T.~Deguchi, N.~Fukushima,
T.~K{\"a}rki, R.~Kashaev, T.~Le, R.~Lickorish, M.~Okamoto, M.~Sakuma and
A.~Vaintrob for their helpful comments.
Thanks are also due to the subscribers to the mailing list `knot'
(http://w3.to/oto/) run by M.~Ozawa.
Finally we thank Maple V (a product of Waterloo Maple Inc.), without the help
of this software the work could have never been carried out.
In fact we used it to check Proposition~\ref{prop:Kashaev=Jones} up to $N=12$
and to find $p$ in Theorem~\ref{thm:Alexander=Jones}.
We also used it for the step-by-step confirmation of the proof of
Proposition~\ref{prop:Kashaev=Jones}.
\end{ack}
\section{Preliminaries}\label{sec:pre}
In this section we recall the definitions of Yang--Baxter operators and
associated link invariants.
If an invertible linear map $R:\C^N\otimes\C^N\to\C^N\otimes\C^N$ satisfies the
following Yang--Baxter equation, it is called a Yang--Baxter operator.
\begin{equation*}
  (R\otimes{id})({id}\otimes{R})(R\otimes{id})
  =
  ({id}\otimes{R})(R\otimes{id})({id}\otimes{R}),
\end{equation*}
where $id:\C^N\to\C^N$ is the identity.
If there exists a homomorphism $\mu:\C^N\to\C^N$ and scalars $\alpha,\beta$
satisfying the following two equations, the quadruple $S=(R,\mu,\alpha,\beta)$
is called an enhanced Yang-Baxer operator \cite{Turaev:INVEM88}.
\begin{gather*}
  (\mu\otimes\mu)R=R(\mu\otimes\mu),
  \\
  \Sp_2\left(R^{\pm1}({id}\otimes\mu)\right)=\alpha^{\pm1}\beta\,{id},
\end{gather*}
where
$\Sp_k:\End(\C^{\otimes{k}})\to\End(\C^{\otimes{k-1}})$ is the
operator trace defined as
\begin{multline*}
  \Sp_k(f)(v_{i_1}\otimes{v_{i_2}}\otimes\dots\otimes{v_{i_{k-1}}})
  \\=
  \sum_{j_1,j_2,\dots,j_{k-1},j=0}^{N-1}
  f_{i_1,i_2,\dots,i_{k-1},j}^{j_1,j_2,\dots,j_{k-1},j}
  (v_{j_1}\otimes{v_{j_2}}\otimes\dots\otimes{v_{j_{k-1}}}\otimes{v_{j}}),
\end{multline*}
where
\begin{equation*}
  f(v_{i_1}\otimes{v_{i_2}}\otimes\dots\otimes{v_{i_k}})
  =
  \sum_{j_1,j_2,\dots,j_k=0}^{N-1}
  f_{i_1,i_2,\dots,i_k}^{j_1,j_2,\dots,j_k}
  \left({v_{j_1}}\otimes{v_{j_2}}\otimes\dots\otimes{v_{j_k}}\right)
\end{equation*}
for a basis $\{v_0,v_1,\dots,v_{N-1}\}$ of $\C^N$.
\par
For an enhanced Yang--Baxter operator one can define a link invariant as follows
\cite{Turaev:INVEM88}.
First we represent a given link $L$ as the closure of a braid $\xi$ with
$n$ strings.
Consider the $n$-fold tensor product $\left(\C^N\right)^{\otimes{n}}$ and
associate the homomorphism
$b_R(B):\left(\C^N\right)^{\otimes{n}}\to\left(\C^N\right)^{\otimes{n}}$ by
replacing $\sigma_i^{\pm1}$
(the usual $i$-th generator of the braid group) in $\xi$ with
\begin{equation*}
  \underset{i-1}{\underbrace{{id}\otimes\dots\otimes{id}}}
  \otimes{R^{\pm1}}\otimes
  \underset{N-i-1}{\underbrace{{id}\otimes\dots\otimes{id}}}.
\end{equation*}
Then taking the operator trace $n$ times we define
\begin{equation*}
  T_S(\xi)
  =
  \alpha^{-w(\xi)}\beta^{-n}
  \Sp_1\left(\Sp_2
         \left(\cdots
           \left(\Sp_n
             \left(b_R(\xi)\mu^{\otimes{n}}
             \right)
           \right)
         \right)
       \right),
\end{equation*}
where $w(\xi)$ is the sum of the exponents.
Then $T_S(\xi)$ defines a link invariant and denoted by $T_S(L)$.
\par
To define the (generalized) Alexander polynomial from an enhanced Yang--Baxter
operator we have to be more careful, since $T_S$ always vanishes in this case.
If the following homomorphism
\begin{equation*}
  T_{S,1}(\xi)
  =
  \alpha^{-w(\xi)}\beta^{-n}
  \Sp_2\left(\Sp_3
         \left(\cdots
           \left(\Sp_n
             \left(b_S(\xi)({id}\otimes\mu^{\otimes(n-1)})
             \right)
           \right)
         \right)
       \right)
   \in\End(\C^N)
\end{equation*}
is a scalar multiple and
\begin{equation*}
  \Sp_1(\mu)T_{S,1}(\xi)=T_S(\xi)
\end{equation*}
for any $\xi$, then the scalar defined by $T_{S,1}(\xi)$ becomes a link
invariant (even if $\Sp_1(\mu)=0$) and is denoted by $T_{S,1}(L)$.
Note that this invariant can be regarded as an invariant for $(1,1)$-tangles,
where a $(1,1)$-tangle is a link minus an open interval.
\section{The intersection of the generalized Alexander polynomials and the
colored Jones polynomials}\label{sec:Alexander}
In \cite{Akutsu/Deguchi/Ohtsuki:JKNOT92} Akutsu, Deguchi and Ohtsuki
defined a generalization of the multivariable Alexander polynomial for colored
links.
First we will briefly describe their construction only for the case where
all the colors are the same according to \cite{Deguchi:1994}.
\par
Fix an integer $N\ge2$ and a complex number $p$.
Put $s=\exp(\pi\sqrt{-1}/N)$ and $[k]=(s^k-s^{-k})/(s-s^{-1})$ for a complex
number $k$.
Note that $[N]=0$ and $[N-k]=[k]$.
\par
Let $U_q(sl(2,\C))$ be the quantum group generated by $X,Y,K$ with the following
relations.
\begin{equation*}
  KX=sXK,\quad KY=s^{-1}YK,\quad XY-YX=\frac{K^2-K^{-2}}{s-s^{-1}}.
\end{equation*}
Let $F(p)$ be the $N$-dimensional vector space over $\C$ with basis
$\{f_0,f_1,\dots,f_{N-1}\}$.
We give an action of $U_q(sl(2,\C))$ on $F(p)$ by
\begin{align*}
  X(f_i)&=\sqrt{[2p-i+1][i]}f_{i-1},
  \\
  Y(f_i)&=\sqrt{[2p-i][i+1]}f_{i+1},
  \\
  K(f_i)&=s^{p-i}f_{i}.
\end{align*}
Using Drinfeld's universal $R$-matrix given in \cite{Drinfeld:ICM86},
we can define a set of enhanced Yang--Baxter operators $S_A(p)$
with complex parameter $p$.
Then Akutsu--Deguchi--Ohtsuki's generalized Alexander
polynomial is defined to be $T_{S_A(p),1}$ by using the notation in the previous
section.
We denote it by $\Phi_N(L,p)$ for a link $L$.
Note that if $N=2$ the invariant $\Phi_2(L,p)$ is the same as the
multivariable Alexander polynomial \cite{Jun:PACJM93}.
\par
Next we review the colored Jones polynomial at the $N$-th root of unity.
There is another $N$-dimensional representation of $U_q(sl(2,\C))$,
corresponding to the usual $N$-dimensional irreducible representation of
$sl(2,\C)$.
Let $E$ be the $N$-dimensional complex vector space with basis
$\{e_0,e_1,\dots,e_{N-1}\}$ and we define the action of $U_q(sl(2,\C))$
by
\begin{align*}
  X(e_i)&=[i+1]e_{i+1},
  \\
  Y(e_i)&=[i]e_{i-1},
  \\
  K(e_i)&=s^{i-(N-1)/2}e_{i}.
\end{align*}
(See for example \cite[(2.8)]{Kirby/Melvin:INVEM91}.)
By using Drinfeld's universal $R$-matrix again we have another enhanced
Yang--Baxter operator $S_J$.
Then the invariant $T_{S_J,1}$ coincides with the colored Jones polynomial
of a link each of whose component decorated by the $N$-dimensional irreducible
representation evaluated at $t=s^2=\exp(2\pi\sqrt{-1}/N)$.
Note that before evaluating at $s^2$, we have to normalize the colored Jones
polynomial so that its value of the trivial knot is one for otherwise
the invariant would be identically zero.
This is well-defined since the colored Jones polynomial defines a well-defined
$(1,1)$-tangle invariant (\cite[(3.9) Lemma]{Kirby/Melvin:INVEM91}).
We will denote $T_{S_J,1}$ by $J_{N}$.
\par
Now we put $p=(N-1)/2$ in the Akutsu--Deguchi--Ohtsuki invariant.
Then since $[N-k]=[k]$, we have
\begin{align*}
  X(f_i)&=[i]f_{i-1},
  \\
  Y(f_i)&=[i+1]f_{i+1},
  \\
  K(f_i)&=s^{(N-1)/2-i}f_{i}
\end{align*}
and so the two representation $F((N-1)/2)$ and $E$ are quite similar.
In fact if we exchange $X$ and $Y$, and replace $K$ with $K^{-1}$ then these two
coincide.
(This automorphism is known as the Cartan automorphism
\cite[p.~123, Lemma VI.1.2]{Kassel:quantum_groups}.)
Therefore they determine the same Yang--Baxter operator and the same
link invariant, that is, we have the following theorem.
\begin{thm}\label{thm:Alexander=Jones}
  The Akutsu--Deguchi--Ohtsuki invariant with all the colors $p=(N-1)/2$
  coincides with the colored Jones polynomial corresponding to the
  $N$-dimensional irreducible representation evaluated at
  $\exp(2\pi\sqrt{-1}/N)$.
  More precisely, we have $\Phi_{N}(L,(N-1)/2)=J_{N}(L)$ for every link $L$.
\end{thm}
\begin{rem}
After finishing this work we were informed by Deguchi that it has already
observed \cite{Deguchi:JPHYS191} that the $R$-matrices given by $F((N-1)/2)$
and $E$ coincide.
\end{rem}
\section{$\cR$-matrix for the colored Jones polynomial at the $N$th root of
         unity}\label{sec:Jones}
Let $R_J$ be the $\cR$-matrix shown in
\cite[Corollary~2.32]{Kirby/Melvin:INVEM91}, which is the $N^2\times N^2$ matrix
with
$((i,j),(k,l))$th entry
\begin{align*}
  \left(R_J\right)_{kl}^{ij}=
  \sum_{n=0}^{\min{(N-1-i,j)}}
  &\delta_{l,i+n}\delta_{k,j-n}
  \frac{(s-s^{-1})^n}{[n]!}
  \frac{[i+n]!}{[i]!}\frac{[N-1+n-j]!}{[N-1-j]!}
  \\
  &\times
  s^{2(i-(N-1)/2)(j-(N-1)/2)-n(i-j)-n(n+1)/2},
\end{align*}
where $[k]!=[k][k-1]\cdots[2][1]$.
Note that our matrix $R_J$ corresponds to $\cR$ in
\cite[Definition~2.35]{Kirby/Melvin:INVEM91}.
This matrix is used to define an enhanced Yang--Baxter operator and
the link invariant $J_{N}$ described in the previous section.
\par
The aim of this section is to transform it to a matrix similar to Kashaev's
$\cR$-matrix.
Let $W$ and $D$ be the $N\times N$ matrices with $(i,j)$th entry
$W_{j}^{i}=s^{2ij}$ and $D_{j}^{i}=\delta_{i,j}s^{(N-1)i}$ respectively,
where $\delta_{i,j}$ is Kronecker's delta.
We will calculate the product
$\tR_J
=(W\otimes W)({id}\otimes{D}){R_J}
({id}\otimes{D^{-1}})(W^{-1}\otimes{W^{-1}})$
with $id$ the $N\times{N}$ identity matrix and show the following proposition.
\begin{prop}\label{prop:J}
\begin{equation*}
  \left(\tR_J\right)_{ab}^{cd}
  =
  \begin{cases}
    \rho(a,b,c,d)(-1)^{a+b+1}\dfrac{[d-c-1]![N-1+c-a]!}{[d-b]![b-a-1]!}
    \quad&\text{if $d\ge b>a\ge c$},
    \\[5mm]
    \rho(a,b,c,d)(-1)^{a+c}\dfrac{[b-d-1]![N-1+c-a]!}{[c-d]![b-a-1]!}
    \quad&\text{if $b>a\ge c\ge d$},
    \\[5mm]
    \rho(a,b,c,d)(-1)^{b+d}\dfrac{[N-1+b-d]![c-a-1]!}{[c-d]![b-a-1]!}
    \quad&\text{if $c\ge d\ge b>a$},
    \\[5mm]
    \rho(a,b,c,d)(-1)^{c+d}\dfrac{[N-1+b-d]![a-b]!}{[c-d]![a-c]!}
    \quad&\text{if $a\ge c\ge d\ge b$},
    \\[5mm]
    0\quad&\text{otherwise},
  \end{cases}
\end{equation*}
where
$\rho(a,b,c,d)=s^{-N^2/2+1/2+c+d-2b+(a-d)(c-b)}[N-1]!(s-s^{-1})^{2(N-1)}/N^2$.
\end{prop}
\begin{proof}
Since
$(W\otimes{W})_{ab}^{ef}=s^{2ae+2bf}$,
$({id}\otimes{D})_{ef}^{kl}=\delta_{e,k}\delta_{f,l}s^{(N-1)l}$,
$\left({id}\otimes{D^{-1}}\right)_{ij}^{gh}
  =\delta_{g,i}\delta_{h,j}s^{-(N-1)j}$,
$\left(W^{-1}\otimes{W^{-1}}\right)_{gh}^{cd}=s^{-2cg-2dh}/N^2$,
we have
\begin{align*}
  &N^2\left(\tR_J\right)_{ab}^{cd}
  \\
  &=
  \sum_{i,j,k,l,e,f,g,h=0}^{N-1}\sum_{n=0}^{\min{(N-1-i,j)}}
  \delta_{e,k}\delta_{f,l}\delta_{g,i}\delta_{h,j}\delta_{l,i+n}\delta_{k,j-n}
  s^{2ae+2bf-2cg-2dh+(N-1)(l-j)}
  \\
  &\quad\times
  \frac{(s-s^{-1})^n}{[n]!}\frac{[i+n]!}{[i]!}\frac{[N-1+n-j]!}{[N-1-j]!}
  s^{2(i-(N-1)/2)(j-(N-1)/2)-n(i-j)-n(n+1)/2}
  \\
  &=
  \sum_{i,j=0}^{N-1}\sum_{n=0}^{\min{(N-1-i,j)}}
  s^{(N-1)^2/2}s^{(2b-2a+N)n-n^2/2-3n/2+(2a-2d+n+2)j+(2b-2c+2j-n)i}
  \\
  &\quad\times
  (s-s^{-1})^n\frac{[N-1+n-j]!}{[N-1-j]!}\qbinom{n+i}{i},
\end{align*}
where
\begin{equation*}
  \qbinom{x}{y}=\frac{[x]!}{[y]![x-y]!}.
\end{equation*}
Since the summation $\sum_{i,j=0}^{N-1}\sum_{n=0}^{\min{(N-1-i,j)}}$ is the same
as $\sum_{n=0}^{N-1}\sum_{j=n}^{N-1}\sum_{i=0}^{N-1-n}$, we have
\begin{align*}
  N^2\left(\tR_J\right)_{ab}^{cd}
  &=
  s^{(N-1)^2/2}\sum_{n=0}^{N-1}\sum_{j=n}^{N-1}
  s^{(2b-2a+N)n-n^2/2-3n/2+(2a-2d+n+2)j}(s-s^{-1})^n
  \\
  &\quad\times
  \frac{[N-1+n-j]!}{[N-1-j]!}S(n,2(b-c+j)-n)
\end{align*}
with $S(\alpha,\beta)=\sum_{i=0}^{N-1}s^{\beta i}\qbinom{\alpha+i}{i}$.
\par
Replacing $j-n$ with $k$, the summation turns out be
$\sum_{k=0}^{N-1}\sum_{n=0}^{N-1-k}$ and we have
\begin{equation*}
  N^2\left(\tR_J\right)_{ab}^{cd}
  =
  s^{(N-1)^2/2}\sum_{k=0}^{N-1}
  s^{2(a-d+1)k}[N-1-k]!X(k),
\end{equation*}
where
\begin{equation*}
  X(k)=
  \sum_{n=0}^{N-1-k}
  (-1)^n s^{2(b-d)n+kn+n(n+1)/2}
  \frac{(s-s^{-1})^n}{[N-1-k-n]!}
  S(n,2(b-c+k)+n).
\end{equation*}
Note that from Lemma~\ref{lem:ST}, we have
\begin{multline}\label{eq:S}
  S(n,2(b-c+k)+n)
  =
  \\
  (1-s^{2(b-c+k-1)})(1-s^{2(b-c+k-2)})\cdots(1-s^{2(b-c+k+n-N+1)}).
\end{multline}
\par
As easily seen from Lemma~\ref{lem:ST}, $S(\alpha,\beta)$ and $T(\alpha,\beta)$
vanishes if $(\beta-\alpha-2)\ge0\ge(\beta+\alpha-2N+2)$ and if
$(\beta+\alpha-1)\ge0\ge(\beta-\alpha+1)$ respectively.
We will use this fact repeatedly from now on and divide the proof into
some cases according to the order of $a,b,c,d$.
\par
First we divide the proof into two cases; $b>c$ and $c\ge b$.
\medskip\par\noindent
{\bf Case 1 ($b>c$).}\qquad
In this case $b-c+k-1\ge 0$ and so from \eqref{eq:S}, $S(n,2(b-c+k)+n)=0$ if
$n\le N-1-k-b+c$.
If $n>N-1-k-b+c$, we see that
\begin{multline*}
  S(n,2(b-c+k)+n)=
  \\
  (-1)^{N-1-n}s^{(N-1-n)(2(b-c+k)+n-N)/2}(s-s^{-1})^{N-1-n}
  \frac{[b-c+k-1]!}{[b-c+k+n-N]!}.
\end{multline*}
Therefore
\begin{align*}
  X(k)&=
  (-1)^{N+1}(s-s^{-1})^{N-1}s^{(b-c+k)(N-1)-N(N-1)/2}
  \frac{[b-c+k-1]!}{[b-c-1]!}
  \\
  &\quad\times
  \sum_{n=N-1-k-b+c}^{N-1-k}
  (-1)^{n}s^{(b+c-2d)n}
  \qbinom{b-c-1}{N-1-k-n}
  \\[5mm]
  &\text{(putting $i=N-1-k-n$)}
  \\
  &=
  (-1)^{k}(s-s^{-1})^{N-1}s^{(2b-2d+k)(N-1)-(b+c-2d)k-N(N-1)/2}
  \frac{[b-c+k-1]!}{[b-c-1]!}
  \\
  &\quad\times
  \sum_{i=0}^{b-c-1}
  (-1)^{i}s^{(2d-b-c)i}
  \qbinom{b-c-1}{i}
  \\
  &=
  (s-s^{-1})^{N-1}s^{(2d-2b)-(b+c-2d+1)k-N(N-1)/2}
  \frac{[b-c+k-1]!}{[b-c-1]!}
  \\
  &\quad\times
  T(b-c-1,2d-b-c),
\end{align*}
with $T(\alpha,\beta)=\sum_{i=0}^{\alpha}(-1)^{i}s^{\beta i}\qbinom{\alpha}{i}$.
From Lemma~\ref{lem:ST} we have
\begin{align*}\label{eq:X}
  X(k)
  &=(s-s^{-1})^{N-1}s^{2(d-b)-(b+c-2d+1)k-N(N-1)/2}
  \frac{[b-c+k-1]!}{[b-c-1]!}
  \\
  &\quad\times
  (1-s^{2(d-c-1)})(1-s^{2(d-c-2)})\cdots(1-s^{2(d-b+1)}).
\end{align*}
We divide the case into two subcases; $d\ge b$ and $b>d$.
\smallskip\par\noindent
{\bf Subcase 1.1 ($d\ge b$).}\quad
Since
\begin{multline*}
  (1-s^{2(d-c-1)})(1-s^{2(d-c-2)})\cdots(1-s^{2(d-b+1)})
  \\
  =
  (-1)^{b-c-1}s^{(2d-b-c)(b-c-1)/2}(s-s^{-1})^{b-c-1}
  \frac{[d-c-1]!}{[d-b]!},
\end{multline*}
we have
\begin{align*}
  X(k)
  &=(-1)^{b+c+1}(s-s^{-1})^{N+b-c-2}
  s^{2(d-b)-N(N-1)/2+(2d-b-c)(b-c-1)/2}
  \\
  &\quad\times
  \frac{[d-c-1]!}{[b-c-1]![d-b]!}
  s^{(2d-b-c-1)k}[b-c+k-1]!
\end{align*}
and so
\begin{align*}
  &N^2\left(\tR_J\right)_{ab}^{cd}
  \\
  &=
  s^{(N^2+1)/2}(-1)^{b+c}(s-s^{-1})^{N+b-c-2}
  s^{2(d-b)-N(N-1)/2+(2d-b-c)(b-c-1)/2}
  \\
  &\quad\times
  \frac{[d-c-1]!}{[d-b]!}
  \sum_{k=0}^{N-1}s^{(2a-b-c+1)k}\frac{[N-1-k]![b-c+k-1]!}{[b-c-1]!}
  \\
  &=
  s^{(N^2+1)/2}(-1)^{b+c}(s-s^{-1})^{N+b-c-2}
  s^{2(d-b)-N(N-1)/2+(2d-b-c)(b-c-1)/2}
  \\
  &\quad\times
  \frac{[d-c-1]![N-1]!}{[d-b]!}
  S(b-c-1,2a-b-c+1)
  \\
  &=
  s^{(N^2+1)/2}(-1)^{b+c}(s-s^{-1})^{N+b-c-2}
  s^{2(d-b)-N(N-1)/2+(2d-b-c)(b-c-1)/2}
  \\
  &\quad\times
  \frac{[d-c-1]![N-1]!}{[d-b]!}
  (1-s^{2(a-b)})(1-s^{2(a-b-1)})\dots(1-s^{2(a-c-N+1)}),
\end{align*}
where the second equality follows from $[N-1-k]!=[N-1]!/[k]!$.
\par
Noting that $a-c-N+1\le 0$, we see that
\begin{multline*}
  (1-s^{2(a-b)})(1-s^{2(a-b-1)})\dots(1-s^{2(a-c-N+1)})
  \\
  =
  \begin{cases}
    0 &\quad\text{if $a\ge b$},
    \\
    (s-s^{-1})^{N+c-b}s^{(2a-b-c-N+1)(N+c-b)/2}\frac{[N-1+c-a]!}{[b-a-1]!}&\quad
    \text{if $b>a$}.
  \end{cases}
\end{multline*}
Therefore we have
\begin{equation*}
  \left(\tR_J\right)_{ab}^{cd}
  =\rho(a,b,c,d)(-1)^{a+b+1}\frac{[d-c-1]![N-1+c-a]!}{[d-b]![b-a-1]!}
\end{equation*}
if $b>a$ and zero otherwise.
Thus the proof is complete for the case where $d\ge b>c$.
(Note that $[N-1+c-a]=0$ if $c>a$.)
\smallskip\par\noindent
{\bf Subcase 1.2 ($b>d$).}\quad
In this case we have
\begin{multline*}
  (1-s^{2(d-c-1)})(1-s^{2(d-c-2)})\cdots(1-s^{2(d-b+1)})
  \\
  =
  s^{(2d-b-c)(b-c-1)/2}(s-s^{-1})^{b-c-1}
  \frac{[b-d-1]!}{[c-d]!}
\end{multline*}
if $c\ge d$ and $0$ otherwise.
Therefore $\left(\tR_J\right)_{ab}^{cd}=0$ if $d>c$.
If $c\ge d$, we have
\begin{align*}
  &N^{2}\left(\tR_J\right)_{ab}^{cd}
  \\
  &\quad=
  s^{(N-1)^2/2}(s-s^{-1})^{N+b-c-2}
  s^{2(d-b)-N(N-1)/2+(2d-b-c)(b-c-1)/2}
  \\
  &\qquad\times
  \frac{[N-1]![b-d-1]!}{[c-d]!}
  S(b-c-1,2a-b-c+1)
  \\
  &\quad=
  s^{-(N-1)/2}(s-s^{-1})^{N+b-c-2}
  s^{2(d-b)+(2d-b-c)(b-c-1)/2}
  \\
  &\qquad\times
  \frac{[N-1]![b-d-1]!}{[c-d]!}
  (1-s^{2(a-b)})(1-s^{2(a-b-1)})\cdots(1-s^{2(a-c-N+1)}).
\end{align*}
This vanishes if $a\ge b$ and equals
$N^2\rho(a,b,c,d)(-1)^{a+c}\frac{[b-d-1]![N-1+c-a]!}{[c-d]![b-a-1]!}$
otherwise, completing the proof for the case $b>d$ and $b>c$ since
$[N-1+c-a]=0$ if $c>a$.
\medskip\par\noindent
{\bf Case 2 ($c\ge b$)}.\quad
First note from \eqref{eq:S}, $X(k)=0$ if $k\ge c-b+1$.
If $k<c-b+1$, we have
\begin{multline*}
  S(n,2(b-c+k)+n)
  \\
  =s^{(2(b-c+k)+n-N)(N-1-n)/2}(s-s^{-1})^{N-1-n}
  \frac{[N-1-n+c-b-k]!}{[c-b-k]!}
\end{multline*}
and so
\begin{align*}
  X(k)
  &=
  (s-s^{-1})^{N-1}s^{(2(b-c+k)-N)(N-1)/2}\frac{[c-b]!}{[c-b-k]!}
  \\
  &\quad\times
  \sum_{n=0}^{N-1-k}s^{(b+c-2d)n}\qbinom{N-1-k-n+c-b}{c-b}
  \\[5mm]
  &\text{(putting $i=N-1-k-n$)}
  \\
  &=
  (s-s^{-1})^{N-1}s^{(2(b-c+k)-N)(N-1)/2+(b+c-2d)(N-1-k)}\frac{[c-b]!}{[c-b-k]!}
  \\
  &\quad\times
  S(c-b,2d-b-c)
  \\
  &=
  (s-s^{-1})^{N-1}s^{(2(b-c+k)-N)(N-1)/2+(b+c-2d)(N-1-k)}\frac{[c-b]!}{[c-b-k]!}
  \\
  &\quad\times
  (1-s^{2(d-c-1)})(1-s^{2(d-c-2)})\cdots(1-s^{2(d-b-N+1)}).
\end{align*}
Since $d-b-N+1\le0$, $X(k)$ vanishes and so does
$\left(\tR_J\right)_{ab}^{cd}$ if $d>c$.
Therefore we assume that $c\ge d$.
\par
In this case since
\begin{multline*}
  (1-s^{2(d-c-1)})(1-s^{2(d-c-2)})\dots(1-s^{2(d-b-N+1)})
  \\
  =s^{(2d-b-c-N)(N-1+b-c)/2}(s-s^{-1})^{N-1+b-c}
  \frac{[N-1+b-d]!}{[c-d]!},
\end{multline*}
we have
\begin{align*}
  X(k)
  &=
  (s-s^{-1})^{2N-2+b-c}
  \\
  &\quad\times
  s^{(2b-2c-N)(N-1)/2+(b+c-2d)(N-1)+(2d-b-c-N)(N-1+b-c)/2}
  \\
  &\quad\times
  \frac{[N-1+b-d]![c-b]!}{[c-d]!}
  \\
  &\quad\times
  s^{(N-1+2d-b-c)k}\frac{1}{[c-b-k]!}.
\end{align*}
Therefore
\begin{align*}
  N^{2}\left(\tR_J\right)_{ab}^{cd}
  &=
  s^{(N-1)^2/2}(s-s^{-1})^{2N-2+b-c}
  \\
  &\quad\times
  s^{(2b-2c-N)(N-1)/2+(b+c-2d)(N-1)+(2d-b-c-N)(N-1+b-c)/2}
  \\
  &\quad\times
  \frac{[N-1]![N-1+b-d]!}{[c-d]!}T(c-b,2a-b-c+1)
  \\
  &=
  s^{(N-1)^2/2}(s-s^{-1})^{2N-2+b-c}
  \\
  &\quad\times
  s^{(2b-2c-N)(N-1)/2+(b+c-2d)(N-1)+(2d-b-c-N)(N-1+b-c)/2}
  \\
  &\quad\times
  (1-s^{2(a-b)})(1-s^{2(a-b-1)})\cdots(1-s^{2(a-c+1)}).
\end{align*}
There are two subcases; $b>a$ and $a\ge b$.
\smallskip\par\noindent
{\bf Subcase 2.1 ($b>a$)}.\quad
Since
\begin{multline*}
  (1-s^{2(a-b)})(1-s^{2(a-b-1)})\cdots(1-s^{2(a-c+1)})
  =
  \\
  s^{(2a-b-c+1)(c-b)}(s-s^{-1})^{c-b}\frac{[c-a-1]!}{[b-a-1]!},
\end{multline*}
we have
\begin{equation*}
  \left(\tR_J\right)_{ab}^{cd}
  =
  \rho(a,b,c,d)(-1)^{b+d}\frac{[N-1+b-d]![c-a-1]!}{[c-d]![b-a-1]!}
\end{equation*}
and the proof for the case where $c\ge b$ and $b>a$ is complete.
(Note that $\left(\tR_J\right)_{ab}^{cd}$ vanishes unless $c\ge d\ge b$.)
\smallskip\par\noindent
{\bf Subcase 2.2 ($a\ge b$)}.\quad
In this case we only have to consider the case where $a\ge c$ for
otherwise $\left(\tR_J\right)_{ab}^{cd}=0$.
Now
\begin{multline*}
  (1-s^{2(a-b)})(1-s^{2(a-b-1)})\cdots(1-s^{2(a-c+1)})
  =
  \\
  (-1)^{b+c}s^{(2a-b-c+1)(c-b)}(s-s^{-1})^{c-b}\frac{[a-b]!}{[a-c]!}
\end{multline*}
and so we have
\begin{equation*}
  \left(\tR_J\right)_{ab}^{cd}
  =
  \rho(a,b,c,d)(-1)^{c+d}\frac{[N-1+b-d]![a-b]!}{[c-d]![a-c]!}.
\end{equation*}
The proof for the case where $c\ge b$ and $a\ge b$ is now complete.
(Note again that $\left(\tR_J\right)_{ab}^{cd}=0$ unless $c\ge d\ge b$.)
\end{proof}
\section{Kashaev's R-matrix and his invariant}\label{sec:Kashaev}
In this section we will calculate Kashaev's $\cR$-matrix given in
\cite{Kashaev:MODPLA95} and prove that it coincides with the matrix $\tR_J$
up to a constant given in the previous section.
\par
We prepare notations following \cite{Kashaev:MODPLA95}.
Fix an integer $N\ge2$.
Put $(x)_n=\prod_{i=1}^{n}(1-x^i)$ for $n\ge 0$.
Define $\theta:\Z\to\{0,1\}$ by
\begin{equation*}
  \theta(n)
  =
  \begin{cases}
    1&\quad\text{if $N>n\ge0$},
    \\
    0&\quad\text{otherwise}.
  \end{cases}
\end{equation*}
For an integer $x$, we denote by $\res(x)\in\{0,1,2,\dots,N-1\}$ the
residue modulo $N$.
\par
Now Kashaev's $\cR$-matrix $R_K$ is given by
\begin{multline*}
  {\left(R_K\right)}_{ab}^{cd}
  \\
  =
  Nq^{1+c-b+(a-d)(c-b)}
  \frac{\theta(\res(b-a-1)+\res(c-d))\theta(\res(a-c)+\res(d-b))}
  {(q)_{\res(b-a-1)}(q^{-1})_{\res(a-c)}(q)_{\res(c-d)}(q^{-1})_{\res(d-b)}}
\end{multline*}
with $q=s^2$.
Note that we are using $P\circ{R}$ with $R$ defined in
\cite[2.12]{Kashaev:MODPLA95} rather than $R$ itself where $P$ is the
homomorphism from $\C^N\otimes\C^N$ to $\C^N\otimes\C^N$ sending
$x\otimes{y}$ to $y\otimes{x}$.
\par
We will show the following proposition.
\begin{prop}\label{prop:K}
\begin{equation*}
  {\left({R}_K\right)}_{ab}^{cd}
    =
    \begin{cases}
    \lambda(a,b,c,d)(-1)^{a+b+1}\dfrac{[d-c-1]![N-1+c-a]!}{[d-b]![b-a-1]!}
    \quad&\text{if $d\ge b>a\ge c$},
    \\[5mm]
    \lambda(a,b,c,d)(-1)^{a+c}\dfrac{[b-d-1]![N-1+c-a]!}{[c-d]![b-a-1]!}
    \quad&\text{if $b>a\ge c\ge d$},
    \\[5mm]
    \lambda(a,b,c,d)(-1)^{b+d}\dfrac{[N-1+b-d]![c-a-1]!}{[c-d]![b-a-1]!}
    \quad&\text{if $c\ge d\ge b>a$},
    \\[5mm]
    \lambda(a,b,c,d)(-1)^{c+d}\dfrac{[N-1+b-d]![a-b]!}{[c-d]![a-c]!}
    \quad&\text{if $a\ge c\ge d\ge b$},
    \\[5mm]
    0\quad&\text{otherwise},
  \end{cases}
\end{equation*}
where
$\lambda(a,b,c,d)=
s^{-N^2/2+N/2+2+c+d-2b+(a-d)(c-b)}(s-s^{-1})^{1-N}N/([N-1]!)^2$.
\end{prop}
\begin{proof}
Since $N-1\ge b-a-1\ge-N$, $N-1\ge c-d\ge -N+1$, $N-1\ge a-c\ge -N+1$ and
$N-1\ge d-b\ge -N+1$, we see that ${({R}_K)}_{ab}^{cd}$ vanishes except for the
following four cases, which have already appeared in Proposition~\ref{prop:J}:
(i) $d\ge b>a\ge c$, (ii) $b>a\ge c\ge d$, (iii) $c\ge d\ge b\ge a$ and
(iv) $a\ge c\ge d\ge b$.
\par
We will only prove the first case because the other cases are similar.
Noting that
\begin{align*}
  (q     )_n&=(-1)^{n}s^{n(n+1)/2}(s-s^{-1})^{n}[n]!,
  \\
  (q^{-1})_n&=s^{-n(n+1)/2}(s-s^{-1})^{n}[n]!,
\end{align*}
we have
\begin{align*}
  {\left({R}_K\right)}_{ab}^{cd}
  &=
  (-1)^{a+b+c+d}Ns^{2+2c-2b+2(a-d)(c-b)+a-d+N/2+1/2}(s-s^{-1})^{-N+1}
  \\
  &\qquad\times
  \frac{1}{[d-b]![b-a-1]![N+c-d]![a-c]!}
\end{align*}
since $\res(b-a-1)=b-a-1$, $\res(a-c)=a-c$, $\res(c-d)=N+c-d$ and
$\res(d-b)=d-b$.
Now since $[N-n]=[n]$, we see that
\begin{equation*}
  \frac{1}{[d-b]![b-a-1]![N+c-d]![a-c]!}
  =
  \frac{1}{([N-1])^2}\frac{[d-c-1]![N+c-a-1]!}{[b-a-1]![d-b]!}.
\end{equation*}
Therefore
\begin{equation*}
  {\left({R}_K\right)}_{ab}^{cd}
  =\lambda(a,b,c,d)(-1)^{a+b+1}\frac{[d-c-1]![N-1+c-a]!}{[d-b]![b-a-1]!}
\end{equation*}
as required.
\end{proof}
\par
Therefore $R_K$ and $R_J$ are equal up to a constant depending
only on $N$.
More precisely we have
\begin{prop}\label{prop:Kashaev=Jones}
Let $R_K$ and $R_J$ be the $\cR$-matrices defined by as above.
Then we have
\begin{equation*}
  {R}_K
  =
  s^{-(N+1)(N-3)/2}(W\otimes{W})({id}\otimes{D})
  R_J({id}\otimes{D^{-1}})({W^{-1}}\otimes{W^{-1}})
\end{equation*}
for any $N\ge2$.
\end{prop}
\begin{proof}
From Propositions~\ref{prop:J} and \ref{prop:K}, we only have to check that
$\rho(a,b,c,d)/\lambda(a,b,c,d)=s^{(N+1)(N-3)/2}$.
We have
\begin{equation*}
  \rho(a,b,c,d)/\lambda(a,b,c,d)
  =(-1)^{N}s^{(N-3)/2}\left(\frac{(s-s^{-1})^{N-1}[N-1]!}{N}\right)^3
\end{equation*}
but this coincides with $s^{(N+1)(N-3)/2}$ as shown below.
\par
We have
\begin{align*}
  (s-s^{-1})^{N-1}[N-1]!
  &=\prod_{k=1}^{N-1}\left(2\sqrt{-1}\sin(k\pi/N)\right)
  \\
  &=\sqrt{-1}^{N-1}\prod_{k=1}^{N-1}\left(2\sin(k\pi/N)\right).
\end{align*}
On the other hand from \cite[I.392-1, p.~33]{Gradshteyn/Ryzhik:1980}, we have
\begin{equation*}
  \sin(Nx)=2^{N-1}\prod_{k=0}^{N-1}\sin(x+k\pi/N).
\end{equation*}
Divided by $\sin x$ and taking the limit $x\to 0$, we have
\begin{equation*}
  N=\prod_{k=1}^{N-1}\left(2\sin(k\pi/N)\right).
\end{equation*}
Therefore we have
\begin{align*}
  (-1)^{N}s^{(N-3)/2}\left(\frac{(s-s^{-1})^{N-1}[N-1]!}{N}\right)^3
  &=
  (-1)^{N}s^{(N-3)/2}\sqrt{-1}^{3(N-1)}
  \\
  &=
  s^{N^2+(N-3)/2+3(N-1)N/2}
  =
  s^{(N+1)(N-3)/2},
\end{align*}
completing the proof.
\end{proof}
\par
We will show that the matrix ${R}_K$ also satisfies the Yang--Baxter
equation.
To do that we prepare a lemma.
\begin{lem}\label{lem:D_through_J}
The matrices $D$ and $D^{-1}$ can go through ${R_J}$ in pair, that is, the
following equality holds.
\begin{equation*}
  ({id}\otimes{D}){R_J}({id}\otimes{D^{-1}})
  =
  (D^{-1}\otimes{id}){R_J}({D}\otimes{id}).
\end{equation*}
\end{lem}
\begin{proof}
It is sufficient to show that $(D\otimes{D})R_J=R_J(D\otimes{D})$.
Since $D_{j}^{i}=\delta_{i,j}s^{(N-1)j}$,
\begin{align*}
  \left((D\otimes{D})R_J\right)_{kl}^{ij}
  &=\sum_{a,b}\delta_{a,k}\delta_{b,l}s^{(N-1)k}s^{(N-1)l}{(R_J)}_{ab}^{ij}
  \\
  &=s^{(N-1)(k+l)}{(R_J)}_{kl}^{ij}
\end{align*}
and
\begin{align*}
  \left(R_J(D\otimes{D})\right)_{kl}^{ij}
  &=\sum_{a,b}\delta_{a,i}\delta_{b,j}s^{(N-1)i}s^{(N-1)j}{(R_J)}_{kl}^{ab}
  \\
  &=s^{(N-1)(i+j)}{(R_J)}_{kl}^{ij}.
\end{align*}
But these two coincide since $(R_J)_{kl}^{ij}$ vanishes unless $k+l=i+j$
(the charge conservation law), completing the proof.
\end{proof}
Using Lemma~\ref{lem:D_through_J} we can give another proof of the following
proposition.
\begin{prop}[Kashaev]\label{prop:YBE_K}
Kashaev's $\cR$-matrix ${R}_K$ satisfies the Yang--Baxter equation, that is,
\begin{equation*}
  ({R}_K\otimes{id})({id}\otimes{{R}_K})({R}_K\otimes{id})
  =
  ({id}\otimes{{R}_K})({R}_K\otimes{id})({id}\otimes{{R}_K}).
\end{equation*}
\end{prop}
\begin{proof}
From Proposition~\ref{prop:Kashaev=Jones}, we have
\begin{align*}
  &(R_K\otimes{id})({id}\otimes{R_K})({R}_K\otimes{id})
  \\
  &\quad=
    (W\otimes{W}\otimes{id})({id}\otimes{D}\otimes{id})(R_J\otimes{id})
    ({id}\otimes{D^{-1}}\otimes{id})(W^{-1}\otimes{W^{-1}}\otimes{id})
  \\
  &\qquad\times
   ({id}\otimes{W}\otimes{W})({id}\otimes{id}\otimes{D})({id}\otimes{R_J})
   ({id}\otimes{id}\otimes{D^{-1}})({id}\otimes{W^{-1}}\otimes{W^{-1}})
  \\
  &\qquad\times
    (W\otimes{W}\otimes{id})({id}\otimes{D}\otimes{id})(R_J\otimes{id})
    ({id}\otimes{D^{-1}}\otimes{id})(W^{-1}\otimes{W^{-1}}\otimes{id})
  \\
  &\quad=
    (W\otimes{W}\otimes{W})({id}\otimes{D}\otimes{D})(R_J\otimes{id})
  \\
  &\qquad\times
    ({id}\otimes{D^{-1}}\otimes{id})({id}\otimes{R_J})
    ({id}\otimes{D}\otimes{id})
  \\
  &\qquad\times
    (R_J\otimes{id})({id}\otimes{D^{-1}}\otimes{D^{-1}})
    (W^{-1}\otimes{W^{-1}}\otimes{W^{-1}})
  \\
  &\quad=
    (W\otimes{W}\otimes{W})({id}\otimes{D}\otimes{D})(R_J\times{id})
  \\
  &\qquad\times
    (id\otimes{id}\otimes{D})({id}\otimes{R_J})
    (id\otimes{id}\otimes{D^{-1}})
  \\
  &\qquad\times
    (R_J\otimes{id})({id}\otimes{D^{-1}}\otimes{D^{-1}})
    (W^{-1}\otimes{W^{-1}}\otimes{W^{-1}})
  \\
  &\quad=
    (W\otimes{W}\otimes{W})(id\otimes{D}\otimes{D^2})
  \\
  &\qquad\times
    (R_J\otimes{id})(id\otimes{R_J})(R_J\otimes{id})
  \\
  &\qquad\times
    (id\otimes{D^{-1}}\otimes{D^{-2}})(W^{-1}\otimes{W^{-1}}\otimes{W^{-1}}).
  \\
  \intertext{Similar calculation shows}
  &({id}\otimes{R_K})(R_K\otimes{id})({id}\otimes{R_K})
  \\
  &\quad=
    (W\otimes{W}\otimes{W})(id\otimes{D}\otimes{D^2})
  \\
  &\qquad\times
    ({id}\otimes{R_J})(R_J\otimes{id})(id\otimes{R_J})
  \\
  &\qquad\times
    (id\otimes{D^{-1}}\otimes{D^{-2}})(W^{-1}\otimes{W^{-1}}\otimes{W^{-1}}).
\end{align*}
From the Yang--Baxter equation for $R_J$ these two coincide, completing
the proof.
\end{proof}
\par
To show that $R_J$ and ${R}_K$ define the same link invariant,
we will construct enhanced Yang--Baxter operators precisely by using them.
\par
Let $\mu_J$ be the $N\times N$-matrix with $(i,j)$-entry
$\left(\mu_J\right)_{j}^{i}=\delta_{i,j}s^{2i-N+1}$.
Then the quadruple $S_J=(R_J,\mu_{J},s^{N^2-1},1)$
is a Yang--Baxter operator and the following lemma holds.
\begin{lem}\label{lem:EYB_J}
\begin{gather*}
  (\mu_J\otimes\mu_J)R_J=R_J(\mu_J\otimes\mu_J),
  \\
  \sum_{j=0}^{N-1}\left((R_J)^{\pm1}(id\otimes\mu_J)\right)_{kj}^{ij}
  =(s^{N^2-1})^{\pm1}id.
\end{gather*}
\end{lem}
\par
Next we will give a Yang--Baxter operator using $R_K$.
Let $\mu_K$ be the $N\times N$-matrix with $(i,j)$-entry
$\left(\mu_K\right)_{j}^{i}=-s\delta_{i,j+1}$.
Then we have
\begin{lem}\label{lem:mu_nu}
\begin{equation*}
  W\,D\,\mu_J\,D^{-1}\,W^{-1}=\mu_K.
\end{equation*}
\end{lem}
\begin{proof}
Since $W_{j}^{a}=s^{2aj}$, $D_{a}^{b}=\delta_{a,b}s^{(N-1)b}$,
$\left(\mu_J\right)_{b}^{c}=\delta_{b,c}s^{2c-N+1}$,
$\left(D^{-1}\right)_{c}^{d}=\delta_{c,d}s^{-(N-1)d}$
and $\left(W^{-1}\right)_{d}^{i}=\delta_{d,i}s^{-2di}/N$, we have
\begin{align*}
  \left(W\,D\,\mu_J\,D^{-1}\,W^{-1}\right)_{j}^{i}
  &=
  \frac{1}{N}(-s)\sum_{a=0}^{N-1}s^{2(j-i+1)a}
  \\
  &=
  -s\delta_{i,j+1},
\end{align*}
completing the proof.
\end{proof}
Combining Lemmas~\ref{lem:EYB_J} and \ref{lem:mu_nu}, we show that
$S_K=({R}_K,\mu_K,-s,1)$ is also an enhanced Yang--Baxter
operator.
\begin{lem}\label{lem:EYB_K}
\begin{gather*}
  (\mu_K\otimes\mu_K){R}_K={R}_K(\mu_K\otimes\mu_K),
  \\
  \sum_{j=0}^{N-1}\left((R_K)^{\pm1}(id\otimes\mu_K)\right)_{kj}^{ij}
  =(-s)^{\pm1}id.
\end{gather*}
\end{lem}
\begin{proof}
Noting that $\mu_J$ and $D$ commutes since they are diagonal, the first equality
follows immediately from that in Lemma~\ref{lem:EYB_J}.
\par
The second equality follows from
\begin{multline*}
  ({R}_K)^{\pm1}({id}\otimes\mu_K)
  \\
  =
  s^{\mp(N+1)(N-3)/2}(W\otimes{W})(id\otimes{D})R_J(id\otimes{\mu_J})
  (id\otimes{D^{-1}})(W^{-1}\otimes{W^{-1}}),
\end{multline*}
completing the proof.
\par
Note that the lemma can also be proved by using
\cite[(2.8) and (2.17)]{Kashaev:MODPLA95}.
\end{proof}
Now we see that $S_J$ and $S_K$ define the same link invariant by using
the following lemma.
\begin{lem}\label{lem:b_K}
Let $\xi$ be an $n$-braid.
Then
\begin{multline*}
  b_{R_K}(\xi)
  =
  \left(W^{\otimes{n}}\right)
  \left(D^{k_1}\otimes{D^{k_2}}\otimes\dots\otimes{D^{k_n}}\right)
  b_{R_J}(\xi)
  \\
  \times
  \left(D^{-k_1}\otimes{D^{-k_2}}\otimes\dots\otimes{D^{-k_n}}\right)
  \left(\left(W^{-1}\right)^{\otimes{n}}\right)
\end{multline*}
for some non-negative integers $k_1,k_2,\dots,k_n$.
\end{lem}
\begin{proof}
In fact we can show that $(k_1,k_2,\dots,k_n)$ is of the form
$$(0,1,2,\dots,d_1,0,1,2,\dots,d_2,\dots,0,1,2,\dots,d_h)$$
by using Lemma~\ref{lem:D_through_J} repeatedly
to `push' $D$ and $D^{-1}$ from left to right
(See the proof of Proposition~\ref{prop:YBE_K}).
Details are omitted.
\end{proof}
\par
Since we know that $J_N=T_{S_J,1}$ is well-defined as described in
\S\ref{sec:pre}, from the previous lemma $T_{S_K,1}(L)$ is also a link
invariant, which we denote by $\langle{L}\rangle_{N}$.
Note that it is implicitly stated in \cite{Kashaev:MODPLA95} that the invariant
can be regarded as an invariant for $(1,1)$-tangles.
Note also that though the invariant was defined only up to a multiple of $s$ in
\cite{Kashaev:MODPLA95}, we can now define it without ambiguity.
\par
Since
\begin{align*}
  &b_{R_K}(\xi)({id}\otimes{\mu_K}^{\otimes{(n-1)}})
  \\
  &\quad=
  \left(W^{\otimes{n}}\right)
  \left(D^{k_1}\otimes{D^{k_2}}\otimes\dots\otimes{D^{k_n}}\right)
  b_{R_J}(\xi)
  \\
  &\qquad\times
  \left({id}\otimes{\mu_J}^{\otimes(n-1)}\right)
  \left(D^{-k_1}\otimes{D^{-k_2}}\otimes\dots\otimes{D^{-k_n}}\right)
  \left(\left(W^{-1}\right)^{\otimes{n}}\right)
\end{align*}
from Lemma~\ref{lem:b_K}
we conclude that $S_J$ and $S_K$ define the same link
invariant.
\begin{thm}
  For any link $L$ and any integer $N\ge2$, $\langle{L}\rangle_{N}$ and
  $J_{N}(L)$ coincide.
\end{thm}
\section{Relation between the simplicial volume and the colored Jones
polynomials.}
Let $K$ be one of the three simplest hyperbolic knots $4_1$, $5_2$ and $6_1$.
Kashaev found in \cite{Kashaev:LETMP97} that the hyperbolic volume of
$S^3\setminus K$, denoted by $\Vol(K)$, coincides numerically with the growth
rate of the absolute value of $\langle K\rangle_{N}$ with respect to $N$.
More precisely,
\begin{equation*}
  \Vol(K)
  =2\pi \, \lim_{N\to \infty}\,
  \frac{\log\left|\langle{K}\rangle_{N}\right|}{N}.
\end{equation*}
\par
We would like to modify his conjecture taking Gromov's simplicial
volume (or Gromov norm) \cite{Gromov:INSHE82} into account.
Let us consider the torus decomposition of the complement of a knot $K$
\cite{Jaco/Shalen:MEMAM79,Johannson:1979}.
Then the simplicial volume of $K$, denoted by $\|K\|$ is equal to the sum
of the hyperbolic volumes of hyperbolic pieces of the decomposition divided
by $v_3$, the volume of the ideal regular tetrahedron in $\bH^3$, the
three-dimensional hyperbolic space.
Recall that it is additive under the connect sum \cite{Soma:INVEM81}
\begin{equation*}
  \|K_1 \sharp K_2\|
  =
  \|K_1\| + \|K_2\|,
\end{equation*}
and that it does not alter by mutation \cite[Theorem~1.5]{Ruberman:INVEM87}.
\par
Noting that $J_{N}$ is multiplicative under the connect sum, that is,
\begin{equation*}
  J_{N}(K_1 \sharp K_2)
  =
  J_{N}(K) \,J_{N}(K_2)
\end{equation*}
and that it does not alter by mutation \cite[Corollary~6.2.5]{Jun:ORSJM89},
we propose the following conjecture.
\begin{conj}[Volume conjecture]
For any knot $K$,
\begin{equation}\label{eq:volume}
  \|K\|
  =
  \frac{2\pi}{v_3}\lim_{N \to \infty}
  \frac{\log\left|J_{N}(K)\right|}{N}.
\end{equation}
\end{conj}
\begin{rem}
First note that if Kashaev's conjecture is true then our conjecture holds
for hyperbolic knots and their connect sums.
It is also true for torus knots since Kashaev and O.~Tirkkonen
[private communication] showed that the right hand side of
\eqref{eq:volume} vanishes in this case by using H.~Morton's formula
\cite{Morton:MATPC95} (see also \cite{Rosso/Jones:JKNOT93}).
\end{rem}
\begin{rem}
Note however that the volume conjecture does not hold for links since $J_{N}$ of
the split union of two links vanishes.
\end{rem}
\par
As a consequence of the volume conjecture, we anticipate the following
simplest case of V.~Vassiliev's conjecture
\cite[6.1 Stabilization conjecture]{Vassiliev:1990}
(see also \cite[Chapter~1, Part V (L), Conjecture]{Kirby:problems}).
\begin{conj}[V.~Vassiliev]\label{conj:Vassiliev}
Assume that every Vassiliev (finite-type) invariant of a knot is
identical to that of the trivial knot.
Then it is unknotted.
\end{conj}
\par
We show that the volume conjecture implies Conjecture~\ref{conj:Vassiliev} by
using the following two lemmas.
\begin{lem}[{\cite[Corollary~4.2]{Gordon:TRAAM83}}]\label{lem:zero_volume}
If $\|K\| = 0$  then $K$ is obtained from the trivial knot
by applying a finite number (possibly zero) of the following two operations:
\begin{enumerate}
\item[1)] making a connect sum,
\item[2)] making a cable.
\end{enumerate}
\end{lem}
\begin{lem}\label{lem:Alexander_of_cable}
If $\|K\| = 0$, then the Alexander polynomial $\Delta(K)$ of $K$ is trivial
if and only if $K$ is the trivial knot.
\end{lem}
\begin{proof}
This lemma comes from Lemma~\ref{lem:zero_volume} and the following three facts
\cite[\S2]{Burau:ABHMS32}
(see also \cite[8.23~Proposition]{Burde/Zieschang:1985}).
\begin{enumerate}
  \item[i)]
    the Alexander polynomial of a non-trivial torus knot is not trivial.
  \item[ii)]
    the Alexander polynomial is multiplicative under the connect sum.
    Therefore if $\Delta(K_1)$ and $\Delta(K_2)$ are non-trivial, then
    $\Delta(K_1\sharp K_2)$  is also non-trivial.
  \item[iii)]
    if $K^\prime$ is a knot obtained from $K$ by a cabling operation, then
    $\Delta(K^\prime)$ is $\Delta(K)f(t)$ with some Laurent polynomial $f(t)$.
    Hence if $\Delta(K)$ is non-trivial, so is $\Delta(K^\prime)$ .
\end{enumerate}
\end{proof}
\par
\begin{proof}
[Proof that the volume conjecture implies Conjecture~\ref{conj:Vassiliev}]
First note that every coefficient of both colored Jones polynomial and Alexander
polynomial as a power series in $h=\log t$ is a Vassiliev invariant.
So a knot $K$ with every Vassiliev invariant trivial has the trivial colored
Jones polynomial for any color and the trivial Alexander polynomial.
In particular $J_{N}(K)=1$ for any $N$.
Therefore assuming the volume conjecture, $\|K\|$ vanishes.
From Lemma~\ref{lem:Alexander_of_cable}, $K$ should be trivial, completing the
proof.
\end{proof}
\begin{rem}
It was pointed out by Vaintrob and Bar-Natan that using the
Melvin--Morton--Rozansky conjecture
\cite{Melvin/Morton:COMMP95,Rozansky:COMMP96} proved by Bar-Natan and
S.~Garoufalidis \cite{BarNatan/Garoufalidis:INVEM96}, we
can also show that a knot is trivial if and only if all of
its colored Jones polynomials are trivial since the Melvin--Morton--Rozansky
conjecture says that the Alexander polynomial can be determined by the colored
Jones polynomials.
\end{rem}
\renewcommand{\thesection}{\Alph{section}}
\setcounter{section}{0}
\section{appendix}
In this appendix we prove some technical formulas used in the paper.
Put
\begin{align*}
  S(\alpha,\beta)&=\sum_{i=0}^{N-1}s^{\beta i}\qbinom{\alpha+i}{i},
  \\
  T(\alpha,\beta)&=\sum_{i=0}^{\alpha}(-1)^{i}s^{\beta i}\qbinom{\alpha}{i}.
\end{align*}
Note that the summation in $S(\alpha,\beta)$ is essentially from $0$ to
$N-1-\alpha$.
Then we have
\begin{lem}\label{lem:ST}
The following formulas hold.
\begin{align*}
  S(\alpha,\beta)&=
    \prod_{j=1}^{N-\alpha-1}(1-s^{\beta-\alpha-2j})
   =(1-s^{\beta-\alpha-2})(1-s^{\beta-\alpha-4})\cdots(1-s^{\beta+\alpha-2N+2}),
  \\
  T(\alpha,\beta)&=
    \prod_{j=1}^{\alpha}(1-s^{\beta+\alpha+1-2j})
    =(1-s^{\beta+\alpha-1})(1-s^{\beta+\alpha-3})\cdots(1-s^{\beta-\alpha+1}).
\end{align*}
\end{lem}
\begin{proof}
We only prove the equality for $S(\alpha,\beta)$ since the other case is
similar.
We use the following quantized Pascal relation.
\begin{equation*}
  \qbinom{\alpha+i}{i}
  =s^{-\alpha}\qbinom{\alpha+i-1}{i-1}+s^{i}\qbinom{\alpha+i-1}{i}.
\end{equation*}
Then since
\begin{align*}
  S(\alpha,\beta)
  &=
  s^{-\alpha}\sum_{i=0}^{N-1}s^{\beta i}\qbinom{\alpha+i-1}{i-1}
  +
  \sum_{i=0}^{N-1}s^{(\beta+1)i}\qbinom{\alpha+i-1}{i}
  \\[5mm]
  &\text{(putting $k=i-1$ in the first term)}
  \\
  &=
  s^{\beta-\alpha}\sum_{k=-1}^{N-2}s^{k}\qbinom{\alpha+k}{k}
  +
  S(\alpha-1,\beta+1)
  \\
  &=
  s^{\beta-\alpha}S(\alpha,\beta)+S(\alpha-1,\beta+1),
\end{align*}
we have the following recursive formula.
\begin{equation*}
  S(\alpha-1,\beta+1)=(1-s^{\beta-\alpha})S(\alpha,\beta).
\end{equation*}
Now the required formula follows since $S(N-1,\gamma)=1$ for any integer
$\gamma$.
\end{proof}
\bibliography{mrabbrev,hitoshi}
\bibliographystyle{amsplain}
\end{document}